\documentclass[12pt,a4paper]{amsart}
\usepackage{amssymb}
\usepackage[mathcal]{eucal}
\usepackage{float}

\theoremstyle{plain}
\newtheorem{thm}{Theorem}[section]
\newtheorem{lem}[thm]{Lemma}
\newtheorem{pro}[thm]{Proposition}
\newtheorem{cor}[thm]{Corollary}

\newtheorem{thmABC}{Theorem}

\newtheorem{proABC}[thmABC]{Proposition}

\newtheorem*{prob}{Problem}

\theoremstyle{remark}

\numberwithin{equation}{section}
\numberwithin{table}{section}

\newcommand{\N}{\mathbb{N}}

\newcommand{\F}{\mathbb{F}}

\renewcommand{\phi}{\varphi} 

\DeclareMathOperator{\Aut}{Aut}
\DeclareMathOperator{\Sym}{Sym}
\DeclareMathOperator{\Alt}{Alt}

\DeclareMathOperator{\Mat}{Mat}
\DeclareMathOperator{\vspan}{span}
\DeclareMathOperator{\wtspec}{wspec}
\DeclareMathOperator{\wt}{wt}
\DeclareMathOperator{\cowt}{co-wt}
\DeclareMathOperator{\com}{com}
\DeclareMathOperator{\supp}{supp}
\DeclareMathOperator{\GL}{GL}
\DeclareMathOperator{\AGL}{AGL}
\DeclareMathOperator{\PSL}{PSL}
\DeclareMathOperator{\PGL}{PGL}
\DeclareMathOperator{\PGamL}{{P}\Gamma{L}}

\begin{document}
\title[Automorphism Groups of Codes]{Automorphism Groups of Cyclic Codes} 

\thanks{\textit{Acknowledgements.} Some of the results in this paper
  were suggested by the first author's doctoral thesis,
  Heinrich-Heine-Universit\"at D\"usseldorf, 2007.}

\author{Rolf Bienert} \address{Mathematisches Institut,
  Heinrich-Heine-Universit\"at D\"usseldorf, 40225 D\"usseldorf,
  Germany} 
\email{Rolf.Bienert@t-online.de}

\author{Benjamin Klopsch} \address{Department of Mathematics, Royal
  Holloway, University of London, Egham TW20 0EX, United Kingdom}
\email{Benjamin.Klopsch@rhul.ac.uk}

\begin{abstract}
  In this article we study the automorphism groups of binary cyclic
  codes.  In particular, we provide explicit constructions for codes
  whose automorphism groups can be described as (a) direct products of
  two symmetric groups or (b) iterated wreath products of several
  symmetric groups.  Interestingly, some of the codes we consider also
  arise in the context of regular lattice graphs and permutation
  decoding.
\end{abstract}



\maketitle

\section{Introduction}

In coding theory one frequently establishes in a natural
way a connection between codes and groups.  For instance, a group
acting on a code may provide valuable insights into the structure of
the code, as with the Mathieu groups acting on the Golay codes.

Let $\mathcal{C}$ be a binary linear code of length $N$ over
$\mathbb{F}_2$.  Up to isomorphism, this simply means that
$\mathcal{C}$ is a subspace of the standard vector space
$\mathbb{F}_2^N$ of dimension $N$ over the prime field $\mathbb{F}_2$
of characteristic $2$.  There is a natural action of the symmetric
group $\Sym(N)$ of degree $N$ on $\mathbb{F}_2^N$ by means of
coordinate permutations.  The automorphism group $\Aut(\mathcal{C})$
of $\mathcal{C}$ is the subgroup of $\Sym(N)$ consisting of all
permutations which map the subspace $\mathcal{C}$ into itself.  (In
general, there are several ways of associating an automorphism group
to a linear code, but the distinctions between these variations
disappear in the context of binary linear codes; cf.\
\cite[Section~1.5]{Hu98}.)  It can be shown that every finite group
arises as the automorphism group of a suitable binary linear code;
cf.~\cite{Ph85}.  The question which finite permutation groups, i.e.\
finite groups with a fixed faithful permutation representation, arise
as automorphism groups of binary linear codes is more subtle; a
possible approach to this problem was indicated in \cite{KnSc80}.

Recall that the binary linear code $\mathcal{C}$ of length $N$ is said
to be cyclic if its automorphism group contains a regular cycle of
length $N$.  The class of binary cyclic codes is both of theoretical
and of practical interest, containing well-known families of codes such
as the quadratic residue codes.  It turns out that the class of groups
which occur as automorphism groups of cyclic codes is much more
restricted.  Indeed, one motivating force behind our work is the
natural and fundamental
\begin{prob}
  Determine the class of finite groups which arise as the automorphism
  groups of (binary) cyclic codes.
\end{prob}  
We bracket the word `binary', because it would be equally interesting
to investigate the problem for other ground fields.  Moreover, one can
ask a corresponding question for permutation groups rather than
groups; cf.\ Theorem~\ref{thm:primitive} below.

\medskip

We now give a summary of our results.  First we exhibit two explicit
families of groups which do not arise as automorphism groups of binary
cyclic codes.

\begin{proABC}\label{pro:cyclic_auto}
  The automorphism group of a binary cyclic code is not isomorphic (as
  an abstract group) to a non-trivial cyclic group of odd order.
\end{proABC}

\begin{thmABC}\label{thm:alternating_auto}
  The automorphism group of a binary cyclic code is not isomorphic (as
  an abstract group) to an alternating group $\Alt(n)$ of degree $n
  \in \{3,4,5,6,7\}$ or $n \geq 9$.  The group $\Alt(8)$ occurs as the
  automorphism group of a binary cyclic code of length $15$.
\end{thmABC}
The exceptional appearance of $\Alt(8)$ can be explained by the
isomorphism $\Alt(8) \cong \PSL(4,2)$; cf.\ our remarks following
Theorem~\ref{thm:primitive}.

Extensive computer calculations show that the automorphism groups of
binary cyclic codes can often be described as iterated wreath products
of symmetric groups; see \cite{Bi07} for a systematic account of
automorphism groups of binary cyclic codes up to length $70$.  We
provide an explicit construction of codes with a prescribed
automorphism group of this type.

\begin{thmABC}\label{thm:wreath_auto}
  Let $r \in \mathbb{N}$, and let $n_1, \ldots, n_r \in
  \mathbb{N}_{\geq 3}$ be odd.  Let $G := \Sym(n_1) \wr \ldots \wr
  \Sym(n_r)$ be the iterated wreath product of symmetric groups of
  degrees $n_1, \ldots, n_r$. Then there exists a binary cyclic code
  $\mathcal{C}$ of length $N := n_1 \cdots n_r$ such that
  $\Aut(\mathcal{C}) \cong G$.
\end{thmABC}

It remains an open problem to construct binary cyclic codes such that
the corresponding automorphism groups are iterated wreath products of
symmetric groups of arbitrary degrees.  The computational evidence
suggests that such products occur frequently, but that extra care must
be taken if the product is to involve as factors the symmetric group
of degree $2$.  A few explicit examples are given in
Section~\ref{sec:examples}.

More rarely, one encounters codes whose automorphism group is a direct
product of two symmetric groups.  Again we are able to offer an
explicit construction of binary linear codes $\mathcal{C}_0(a,b)$,
parameterised by $a,b \in \N$ with $a \leq b$, whose automorphism
groups are of this type; for exceptional values of $a,b$ the
automorphism groups are, in fact, slightly larger.  A detailed
description of the family of codes $\mathcal{C}_0(a,b)$ is given in
Propositions~\ref{pro:C_0_dim_md}, \ref{pro:C_0_auto} and
Corollary~\ref{cor:C_0_cyclic}.  As a consequence we record

\begin{thmABC}\label{thm:direct}
  Let $a,b \in \mathbb{N}$ with $2 < a < b$ and $\gcd(a,b) = 1$.  Then
  there exist binary cyclic codes $\mathcal{C}$ such that
  $\Aut(\mathcal{C}) \cong \Sym(a) \times \Sym(b)$.
\end{thmABC}

Interestingly, some of the codes $\mathcal{C}_0(a,b)$ were recently
studied by Key and Seneviratne in the context of regular lattice
graphs and permutation decoding.  In fact, we provide a new, unified
treatment of a related family $\mathcal{C}_1(a,b)$ of binary linear
codes whose study was initiated in~\cite{KeSe07}.  Our approach leads
to a complete description of the automorphism groups of these codes,
allowing us, for instance, to decide which of the codes
$\mathcal{C}_1(a,b)$ are cyclic.  A detailed description of the family
of codes $\mathcal{C}_1(a,b)$ is given in
Propositions~\ref{pro:C_1_equal_C_0}, \ref{pro:C_1_dim_md} and
\ref{pro:C_1_auto}.

Finally, we use results from the well-developed theory of permutation
groups and modular permutation representations to give a description
of the primitive permutation groups which occur as automorphism groups
of binary cyclic codes.

\begin{thmABC}\label{thm:primitive}
  Let $G \leq \Sym(N)$ be the automorphism group of a binary cyclic
  code $\mathcal{C}$, and suppose that $G$ is a primitive permutation
  group.  Then one of the following holds.
  \begin{enumerate}
  \item $C_p \lneqq G \lneqq \AGL(1,p)$ where $p = N \geq 5$ is a
    prime.
  \item $G = \Sym(N)$; in this case $\mathcal{C}$ is one of four
    elementary codes.
  \item $G = \PGamL(d,q)$ where $d \geq 3$, $q = 2^k$ for $k \in \N$
    and $N = (q^d - 1)/(q-1)$.
  \item $G = M_{23}$ and $N = 23$.
  \end{enumerate}
  Moreover, each of the groups listed in \textup{(2)}--\textup{(4)}
  does occur as the automorphism group of a suitable binary cyclic
  code.
\end{thmABC}

It remains an open problem to find out precisely which subgroups of
affine groups occur as automorphism groups of binary cyclic codes.
This appears to be essentially a question in combinatorial number
theory.  Computer calculations show that, for instance, there exists a
binary cyclic $[17,8,6]$-code whose automorphism group is $C_8 \ltimes
C_{17} \leq \AGL(1,17)$.  More explicit examples are given in
Section~\ref{sec:examples}.

Part (3) of Theorem~\ref{thm:primitive} can be regarded as a
generalisation of the well-known fact that the automorphism group of
the binary Hamming code of length $2^d - 1$ is $\PSL(d,2)$.

\medskip

\noindent
\textit{Organisation.} The paper is divided into eight sections.
Section~\ref{sec:preliminary} contains a brief summary of general
notions and terminology as well as the constructions of the specific
code families $\mathcal{C}_0(a,b)$, $\mathcal{C}_1(a,b)$ and
$\mathcal{K}(n_1,\ldots,n_r)$.  Propositions~\ref{pro:C_0_dim_md},
\ref{pro:C_0_auto} and Corollary~\ref{cor:C_0_cyclic} in
Section~\ref{sec:C_0} describe the structure of the codes
$\mathcal{C}_0(a,b)$ and imply Theorem~\ref{thm:direct}.
Propositions~\ref{pro:C_1_equal_C_0}, \ref{pro:C_1_dim_md} and
\ref{pro:C_1_auto} in Section~\ref{sec:C_1} describe the structure of
the codes $\mathcal{C}_1(a,b)$ related to rectangular lattice graphs.
In Section~\ref{sec:wreath} we determine the automorphism groups of
the codes $\mathcal{K}(n_1,\ldots,n_r)$ and thereby prove
Theorem~\ref{thm:wreath_auto}.  Proposition~\ref{pro:cyclic_auto} and
Theorem~\ref{thm:alternating_auto} are established in
Section~\ref{sec:cyclic_alt}.  In Section~\ref{sec:primitive} we
describe the primitive permutation groups which occur as automorphism
groups of binary cyclic codes, thus proving
Theorem~\ref{thm:primitive}.  Finally, in Section~\ref{sec:examples}
we give several examples of binary cyclic codes with automorphism
groups which are not fully explained by the results in this paper.

%


\section{Preliminaries and basic set-up}\label{sec:preliminary}

\subsection{General notions} \label{sec:prel_general} Let $\Omega$ be
a finite set of size $N := \lvert \Omega \rvert$.  Consider an
$N$-dimensional vector space $\mathcal{V}$ over the field $\F_2$, with
a fixed \emph{standard basis} $\mathbf{e}_\omega$ indexed by $\omega
\in \Omega$.  Binary linear codes $\mathcal{C}$ of length $N$ can then
be constructed as subspaces of $\mathcal{V}$ with respect to the
standard basis.  As the standard basis is indexed by elements of the
set $\Omega$ we regard the automorphism group $\Aut(\mathcal{C})$ of
any linear code $\mathcal{C} \leq \mathcal{V}$ as a subgroup of
$\Sym(\Omega) \cong \Sym(N)$.

The \emph{support} and the \emph{weight} of $\mathbf{v} = \sum_{\omega
  \in \Omega} v_\omega \mathbf{e}_\omega \in \mathcal{V}$ are defined
as
$$
\supp(\mathbf{v}) := \{ \omega \in \Omega \mid v_\omega \not = 0 \}
\quad \text{and} \quad \wt(\mathbf{v}) := \lvert \supp(\mathbf{v})
\rvert.
$$
The \emph{common weight} of $\mathbf{v},\mathbf{w} \in \mathcal{V}$ is
defined as
$$
\com(\mathbf{v},\mathbf{w}) := \lvert \supp(\mathbf{v}) \cap
\supp(\mathbf{w}) \rvert,
$$  
Clearly, the weight and common weight functions are invariant under
coordinate permutations.  The \emph{weight spectrum} and the
\emph{minimum distance} of a linear code $\mathcal{C} \leq
\mathcal{V}$ are given by
$$
\wtspec(\mathcal{C}) := \{ \wt(\mathbf{v}) \mid \mathbf{v} \in
\mathcal{C} \} \quad \text{and} \quad d(\mathcal{C}) := \min \left(
\wtspec(\mathcal{C}) \setminus \{0\} \right).
$$
The \emph{co-weight} of $\mathbf{v} \in \mathcal{V}$ is defined as
$$
\cowt(\mathbf{v}) := \lvert \Omega \setminus \supp(\mathbf{v}) \rvert
= N - \wt(\mathbf{v}).
$$
We call $\widehat d (\mathcal{C}) := \min \{ \cowt(\mathbf{v}) \mid
\mathbf{v} \in \mathcal{C} \}$ the \emph{minimum co-distance} of
$\mathcal{C}$.

We call an element $\mathbf{v} \in \mathcal{C}$ \emph{decomposable} if
it can be written as $\mathbf{v} = \mathbf{w}_1 + \mathbf{w}_2$ where
$\mathbf{w}_1, \mathbf{w}_2 \in \mathcal{C} \setminus \{ \mathbf{0}
\}$ with $\supp(\mathbf{w}_1) \cap \supp(\mathbf{w}_2) = \varnothing$.
An element of $\mathcal{C}$ is \emph{indecomposable} if it is non-zero
and not decomposable.  Clearly, any element of $\mathcal{C}$ of minimum weight
$d(\mathcal{C})$ is indecomposable, and the set of indecomposable
elements is invariant under $\Aut(\mathcal{C})$.  

\subsection{The codes $\mathcal{C}_0(a,b)$ and $\mathcal{C}_1(a,b)$}
\label{sec:prel_C0_C1} Let $a,b \in \N$ with $a \leq b$, and set
$\Omega := \{1,\ldots,a\} \times \{1,\ldots,b\}$ so that $N := \lvert
\Omega \rvert = ab$.  Consider the $N$-dimensional vector space
$\mathcal{V} := \Mat(a,b,\F_2)$ of all $a \times b$ matrices over the
field $\F_2$.  As a \emph{standard basis} of $\mathcal{V}$ we fix
$$
\{ \mathbf{e}_{ij} \mid (i,j) \in \Omega \}, \qquad \text{where
  $\mathbf{e}_{ij} := (\delta_{ik}\delta_{jl})_{kl} \in
  \Mat(a,b,\F_2)$}
$$
denotes the elementary matrix whose $(k,l)$-entry equals $1$ if
$(k,l)=(i,j)$ and $0$ otherwise.

In order to construct specific binary linear codes
$\mathcal{C}_0(a,b)$ and $\mathcal{C}_1(a,b)$ of length $N$ as
subspaces of $\mathcal{V}$ we define the \emph{elementary row
  matrices}
\begin{align*}
  \mathbf{r}_i & := \sum_{j=1}^b \mathbf{e}_{ij} = %
  \bordermatrix{ & & & & & & \cr %
    & 0 & 0 & 0 & \cdots & \cdots & 0 \cr %
    {\scriptstyle i \rightarrow} & 1 & 1 & 1 & \cdots & \cdots & 1
    \cr %
    & 0 & 0 & 0 & \cdots & \cdots & 0 } %
  \qquad \text{for $i\in \{1,\ldots,a\}$,} \\
  \intertext{and the \emph{elementary column matrices}} \mathbf{c}_j &
  := \sum_{i=1}^a \mathbf{e}_{ij} = %
  \bordermatrix{ & & & & {\scriptstyle j \downarrow} & & & \cr %
    & 0 & \cdots & 0 & 1 & 0 & \cdots & 0 \cr %
    & \vdots & & \vdots & \vdots & \vdots & & \vdots \cr %
    & 0 & \cdots & 0 & 1 & 0 & \cdots & 0 } %
  \qquad \text{for $j \in \{1,\ldots,b\}$.}
\end{align*}
Writing $\mathbf{R} := \{ \mathbf{r}_i \mid 1 \leq i \leq a\}$ and
$\mathbf{C} := \{ \mathbf{c}_j \mid 1 \leq j \leq b\}$, we define
$$
\mathcal{C}_0 := \mathcal{C}_0(a,b) := \vspan \langle \mathbf{R} \cup
\mathbf{C} \rangle
$$
to be the vector subspace of $\mathcal{V}$ spanned by the elementary
row and column matrices.  Basic invariants of the code $\mathcal{C}_0$
and the structure of its automorphism group are determined in
Section~\ref{sec:C_0}.  Here we record an inherent symmetry in the
construction of $\mathcal{C}_0$: we notice that $\Aut(\mathcal{C}_0)$
contains the group $\Sym(a) \times \Sym(b)$ which embeds into
$\Sym(\Omega)$ via the imprimitive action
\begin{equation}\label{equ:action}
(i,j)^{(\sigma,\tau)} = (i^\sigma,j^\tau) \qquad \text{for $(i,j) \in
  \Omega$ and $(\sigma,\tau) \in \Sym(a) \times \Sym(b)$.}
\end{equation}
Indeed, in the corresponding action on $\mathcal{V}$, the first factor
$\Sym(a)$ permutes the elements of $\mathbf{R}$ among themselves and
fixes each elementary column matrix, whereas the second factor
$\Sym(b)$ permutes the elements of $\mathbf{C}$ and fixes each
elementary row matrix.

Interestingly, the code $\mathcal{C}_0 = \mathcal{C}_0(a,b)$ also
arises naturally in the study of binary codes defined from rectangular
lattice graphs by Key and Seneviratne~\cite{KeSe07}.  They associate a
binary code $\mathcal{C}_1 = \mathcal{C}_1(a,b)$ to the line graph
$L_2(a,b)$ of the complete bipartite graph $K_{a,b}$ and show that for
such a code permutation decoding can be used for full
error-correction.  A key observation is that $\Aut(\mathcal{C}_1)$
contains $\Sym(a) \times \Sym(b)$.  Using the notation introduced
above, one easily checks that
$$
\mathcal{C}_1 = \mathcal{C}_1(a,b) = \vspan \langle \mathbf{r}_i +
\mathbf{c}_j \mid (i,j) \in \Omega \rangle
$$ 
from which the inclusion $\Sym(a) \times \Sym(b) \subseteq
\Aut(\mathcal{C}_1)$ is now obvious.

\subsection{A weight formula for elements of $\mathcal{C}_0(a,b)$}
\label{sec:weight-formula}

For later use we record a weight formula for elements of
$\mathcal{C}_0$ and the weight spectra of $\mathcal{C}_0$,
$\mathcal{C}_1$.  Let $\mathbf{v} \in \mathcal{C}_0$.  Then there are
$X \subseteq \{1,\ldots,a\}$ and $Y \subseteq \{1,\ldots,b\}$ such
that $\mathbf{v} = \sum_{i \in X} \mathbf{r}_i + \sum_{j \in Y}
\mathbf{c}_j$.  Writing $x := \lvert X \rvert$ and $y := \lvert Y
\rvert$, we find $\sigma \in \Sym(a)$ and $\tau \in \Sym(b)$ such that
$X^\sigma = \{1,\ldots,x\}$ and $Y^\tau = \{1,\ldots,y\}$.  Since the
weight function is invariant under the action of $\Sym(a) \times
\Sym(b)$ on $\mathcal{V}$, corresponding to the action on $\Omega$
described in \eqref{equ:action}, this yields
\begin{equation}\label{equ:weight-formula}
\wt(\mathbf{v}) = \wt(\mathbf{v}^{(\sigma,\tau)}) = \wt \left(
  \sum_{i=1}^x \mathbf{r}_i + \sum_{j=1}^y \mathbf{c}_j \right) =
(a-x)y + (b-y)x.
\end{equation}
Observe that $\mathbf{v} \in \mathcal{C}_1$ if and only if $x+y
\equiv_2 0$.  Thus we obtain
\begin{align*}
  \wtspec(\mathcal{C}_0) &= \{ (a-x)y + (b-y)x \mid 0 \leq x \leq a, \,
  0 \leq y \leq b \}, \\
  \wtspec(\mathcal{C}_1) &= \{ (a-x)y + (b-y)x \mid 0 \leq x \leq a, \,
  0 \leq y \leq b, \, x+y \equiv_2 0 \}.
\end{align*}

\subsection{The codes $\mathcal{K}(n_1,\ldots,n_r)$ }
\label{sec:prel_K} Let $r \in \N$, and let $n_1,\ldots,n_r \in
\N_{\geq 3}$ be odd.  We put $n_0 := 1$.  In this subsection we
provide a recursive definition for a sequence of codes $\mathcal{K}_i
= \mathcal{K}(n_1,\ldots,n_i)$, $i \in \{1,\ldots,r\}$, whose
automorphism groups are later shown to be iterated wreath products of
symmetric groups; see Theorem~\ref{thm:auto_K}.  Set
$$
\Omega_0 := \{1\} \quad \text{and} \quad \Omega_i := \{1,\ldots,n_i\}
\times \Omega_{i-1} \text{ for $i \in \{1,\ldots,r\}$}.
$$
For $i \in \{0,\ldots,r\}$ we fix an $\F_2$-vector space
$$
\mathcal{V}_i := \bigoplus_{\omega \in \Omega_i} \F_2
\mathbf{e}_\omega, \qquad \dim(\mathcal{V}_i) = \lvert
\Omega_i \rvert = n_1 \cdots n_i,
$$
with standard basis $\mathbf{e}_\omega$ indexed by $\omega \in
\Omega_i$ and we set
$$
\mathcal{A}_i := \vspan \langle \mathbf{a}_i^{(k,l)} \mid k,l \in
\{1,\ldots,n_i\} \text{ with } k \not = l \rangle \leq \mathcal{V}_i
$$
where $\mathbf{a}_i^{(k,l)} := \sum_{\omega \in \Omega_{i-1}}
\mathbf{e}_{(k,\omega)} + \mathbf{e}_{(l,\omega)}$ for any distinct
$k,l \in \{1,\ldots,n_i\}$.  We put $\mathcal{K}_0 := \{0\}$, and for
$i \in \{1,\ldots,r\}$ we define recursively
\begin{equation}\label{equ:defK}
  \mathcal{K}_i := \mathcal{K}(n_1,\ldots,n_i) :=
  \begin{cases}
    \mathcal{K}_{i-1}^{(1)} \oplus \ldots \oplus
    \mathcal{K}_{i-1}^{(n_i)} \oplus \mathcal{A}_i & \text{if $i
      \equiv_2 1$,} \\
    \mathcal{K}_{i-1}^{(1)} \oplus \ldots \oplus
    \mathcal{K}_{i-1}^{(n_i)} & \text{if $i \equiv_2 0$,}
  \end{cases}
\end{equation}
where for each $k \in \{1,\ldots,n_i\}$ the summand
\begin{equation}\label{equ:copy}
\mathcal{K}_{i-1}^{(k)} := \left\{ \sum\nolimits_{\omega \in \Omega_{i-1}}
  c_\omega \mathbf{e}_{(k,\omega)} \mid \sum\nolimits_{\omega \in\Omega_{i-1}}
  c_\omega \mathbf{e}_\omega \in \mathcal{K}_{i-1} \right\} \leq
\mathcal{V}_i
\end{equation}
is an isomorphic copy of $\mathcal{K}_{i-1}$ 
with 
$$
\supp(\mathcal{K}_{i-1}^{(k)}) = \{ k \} \times \Omega_{i-1} \qquad
\text{if $i \geq 2$.}
$$
The directness of the sums in \eqref{equ:defK} will be justified in
the proof of the following proposition.

\begin{pro} \label{pro:basicsK}
  Let $i \in \{1,\ldots,r\}$.  Then $\mathcal{K}_i =
  \mathcal{K}(n_1,\ldots,n_i)$ is a binary cyclic code of length $n_1
  \cdots n_i$. It has minimum distance $d(\mathcal{K}_i) = 2$ and
  minimum co-distance $\widehat d(\mathcal{K}_i) = \prod_{j=1}^{\lfloor
    i/2 \rfloor} n_{2j}$.  Its dimension is
  $$
  \dim(\mathcal{K}_i) = 
  \begin{cases}
    \sum_{j=1}^{i+1} (-1)^{j+1} \prod_{k=j}^i n_k & \text{if $i
      \equiv_2
      1$,} \\
    \sum_{j=1}^i (-1)^{j+1} \prod_{k=j}^i n_k & \text{if $i \equiv_2
      0$.}
  \end{cases}
  $$
\end{pro}

\begin{proof}
  Clearly, $\mathcal{K}_i$ is a binary linear code of length $\lvert
  \Omega_i \rvert = n_1 \cdots n_i$.  A short induction argument shows
  that $\Aut(\mathcal{K}_i)$ contains the iterated wreath product
  $\Sym(n_1) \wr \ldots \wr \Sym(n_i)$, in its natural imprimitive
  action on $\Omega_i$.  The wreath product contains a regular cyclic
  subgroup; cf.\ the treatment of the case $2=a<b$ in the proof of
  Corollary~\ref{cor:C_0_cyclic}.  Hence $\mathcal{K}_i$ is a cyclic
  code.

  A straightforward induction shows that
  $$
  d(\mathcal{K}_i) = 2 \quad \text{and} \quad \mathcal{K}_i \subseteq
  \{ \mathbf{v} \in \mathcal{V}_i \mid \wt(\mathbf{v}) \equiv_2 0 \}.
  $$
  Next we comment on the directness of the sums in \eqref{equ:defK}.
  The sum $\mathcal{B}_i := \mathcal{K}_{i-1}^{(1)} \oplus \ldots
  \oplus \mathcal{K}_{i-1}^{(n_i)}$ is direct, since
  $\supp(\mathcal{K}_{i-1}^{(k)})$ and
  $\supp(\mathcal{K}_{i-1}^{(l)})$ are disjoint for any distinct $k,l
  \in \{1,\ldots,n_i\}$.  It remains to explain that $\mathcal{B}_i
  \cap \mathcal{A}_i = \{\mathbf{0}\}$.  Observe that if $\mathbf{a}
  \in \mathcal{A}_i$ is non-zero, then there exists $k \in
  \{1,\ldots,n_i\}$ such that $\supp(\mathcal{K}_{i-1}^{(k)})
  \subseteq \{k\} \times \Omega_{i-1} \subseteq \supp(\mathbf{a})$.
  Hence $\mathbf{a} \in \mathcal{B}_i$ would imply $\sum_{\omega \in
    \Omega_{i-1}} \mathbf{e}_{(k,\omega)} \in
  \mathcal{K}_{i-1}^{(k)}$.  But $\wt(\sum_{\omega \in \Omega_{i-1}}
  \mathbf{e}_{(k,\omega)}) = \lvert \Omega_{i-1} \rvert \equiv_2 1$,
  whereas $\wt(\mathbf{v}) \equiv_2 0$ for any $\mathbf{v} \in
  \mathcal{K}_{i-1}^{(k)}$.  Therefore $\mathbf{a} \not \in
  \mathcal{B}_i$, and $\mathcal{B}_i \cap \mathcal{A}_i =
  \{\mathbf{0}\}$.

  From this we can easily compute the dimension of $\mathcal{K}_i$.
  We have $\dim(\mathcal{K}_i) = n_i \dim(\mathcal{K}_{i-1}) +
  \dim(\mathcal{A}_i) = n_i \dim(\mathcal{K}_{i-1}) + (n_i -1)$ if $i
  \equiv_2 1$, and $\dim(\mathcal{K}_i) = n_i \dim(\mathcal{K}_{i-1})$
  if $i \equiv_2 0$.  Induction gives the desired formula.

  Finally, we determine the minimum co-distance of $\mathcal{K}_i$.
  We contend that $\widehat d(\mathcal{K}_i) = \widehat d(\mathcal{K}_{i-1})$
  if $i \equiv_2 1$, and $\widehat d(\mathcal{K}_i) = n_i \widehat
  d(\mathcal{K}_{i-1})$ if $i \equiv_2 0$.  Induction then yields the
  desired formula.  For $i \equiv_2 0$ our claim follows directly from
  \eqref{equ:defK}.  Now suppose that $i \equiv_2 1$.  Recalling that
  $n_i \equiv_2 1$, it is not difficult to see that a typical element
  realising minimum co-distance in $\mathcal{K}_i$ is $\mathbf{v} =
  \mathbf{v}^{(1)} + \sum_{j=1}^{\lfloor n_i /2 \rfloor}
  \mathbf{a}_i^{(2j,2j+1)}$ where $\mathbf{v}^{(1)} \in
  \mathcal{K}_{i-1}^{(1)}$ corresponds to an element realising minimum
  co-distance in $\mathcal{K}_{i-1}$.
\end{proof}


\section{The codes $\mathcal{C}_0(a,b)$ and their automorphism groups}
\label{sec:C_0}

Let $a,b \in \N$ with $a \leq b$.  We make free use of the notation
introduced in Sections~\ref{sec:prel_general}, \ref{sec:prel_C0_C1}
and \ref{sec:weight-formula}.  The aim of this section is to establish
the following results concerning the binary linear code $\mathcal{C}_0
= \mathcal{C}_0(a,b)$ and its automorphism group.

\begin{pro}\label{pro:C_0_dim_md}
  The code $\mathcal{C}_0 = \mathcal{C}_0(a,b)$ has dimension
  $\dim(\mathcal{C}_0) = a+b-1$ and minimum distance $d(\mathcal{C}_0)
  = a$.

  The special case $a \in \{1,2\}$ allows the following explicit
  description.
  \begin{enumerate}
  \item If $1=a \leq b$, then $\mathcal{C}_0 =\mathcal{V}$.
  \item If $2=a=b$, then $\mathcal{C}_0 = \{ \mathbf{v} \in
    \mathcal{V} \mid \wt(\mathbf{v}) \equiv_2 0 \}$.
  \item If $2=a<b$, then $\mathcal{C}_0 = \{ \sum c_{ij}
    \mathbf{e}_{ij} \mid \forall j,k : c_{1j} + c_{2j} = c_{1k} +
    c_{2k} \}$.
  \end{enumerate}
\end{pro}

\begin{pro}\label{pro:C_0_auto}
  Let $\mathcal{C}_0 = \mathcal{C}_0(a,b)$ as above.
  \begin{enumerate}
  \item If $1=a \leq b$, then $\Aut(\mathcal{C}_0) = \Sym(\Omega)
    \cong \Sym(b)$.
  \item If $2=a=b$, then $\Aut(\mathcal{C}_0) = \Sym(\Omega) \cong
    \Sym(4)$.
  \item If $2=a<b$, then $\Aut(\mathcal{C}_0) = C_2 \wr \Sym(b)$.
  \item If $2<a=b$, then $\Aut(\mathcal{C}_0) = \Sym(a) \wr
    C_2$.
  \item If $2<a<b$, then $\Aut(\mathcal{C}_0) = \Sym(a) \times
    \Sym(b)$.
  \end{enumerate}
\end{pro}

As a corollary, we record for which values of $(a,b)$ the code
$\mathcal{C}_0$ is cyclic, i.e.\ for which $(a,b)$ the permutation
group $\Aut(\mathcal{C}_0) \leq \Sym(\Omega)$ contains a regular
cyclic subgroup.

\begin{cor}\label{cor:C_0_cyclic}
  Let $\mathcal{C}_0 = \mathcal{C}_0(a,b)$ as above.  Then
  $\mathcal{C}_0$ is cyclic, if and only if $a \in \{1,2\}$ or
  $\gcd(a,b) = 1$.
\end{cor}

Note that Theorem~\ref{thm:direct} follows from
Proposition~\ref{pro:C_0_auto} and Corollary~\ref{cor:C_0_cyclic}.  We
now supply the proofs of the stated results.

\begin{proof}[Proof of Proposition~\ref{pro:C_0_dim_md}]
  First we determine the dimension of $\mathcal{C}_0 = \vspan \langle
  \mathbf{R} \cup \mathbf{C} \rangle$.  For each $(i,j) \in \Omega$
  there are precisely two elements in $\mathbf{R} \cup \mathbf{C}$
  which have a non-zero entry in the $(i,j)$-position, namely
  $\mathbf{r}_i$ and $\mathbf{c}_j$.  Therefore $\sum_{i=1}^a
  \mathbf{r}_i + \sum_{j=1}^b \mathbf{c}_j = \mathbf{0}$ and this is
  the only non-trivial linear dependence relation among the $a+b$
  elementary row and column matrices.  Hence $\dim \mathcal{C}_0 =
  a+b-1$.

  In order to determine $d(\mathcal{C}_0)$ we employ the weight
  formula \eqref{equ:weight-formula}.  Let $\mathbf{v} = \sum_{i \in
    X} \mathbf{r}_i + \sum_{j \in Y} \mathbf{c}_j$ and write $x :=
  \lvert X \rvert$, $y := \lvert Y \rvert$, as in
  Subsection~\ref{sec:weight-formula}.  If $x \in \{1,\ldots,a-1\}$,
  then 
  $$
  \wt(\mathbf{v}) = (a-x)y + (b-y)x \geq y + (b-y) = b,
  $$ with equality if and only if $a=2$ or $(x,y) \in \{ (1,0),(a-1,b)
  \}$.  In the latter case, $\mathbf{v} \in \mathbf{R}$.  If $x \in
  \{0,a\}$, then $\wt(\mathbf{v}) = (a-x)y + (b-y)x$ is a multiple of
  $a$, and equal to $a$ if and only if $(x,y) \in \{ (0,1),(a,b-1)\}$,
  equivalently $\mathbf{v} \in \mathbf{C}$.  This analysis shows, in
  particular, that $d(\mathcal{C}_0) = a$.

  It remains to justify the explicit description of $\mathcal{C}_0$
  for $a \in \{1,2\}$.  The cases $1=a \leq b$ and $2=a=b$ are easily
  dealt with, noting that $\dim(\mathcal{C}_0) = b = N$ and
  $\dim(\mathcal{C}_0) = 3 = N-1$, respectively.  Now suppose that
  $2=a<b$.  The description of $\mathcal{C}_0$ can be checked by
  counting: clearly, $\{ \sum c_{ij} \mathbf{e}_{ij} \mid \forall j,k
  : c_{1j} + c_{2j} = c_{1k} + c_{2k} \} \subseteq \mathcal{C}_0$ and
  both sets contain the same number of elements, namely $2^b + 2^b =
  2^{a+b-1} = 2^{\dim(\mathcal{C}_0)}$.
\end{proof}

\begin{proof}[Proof of Proposition~\ref{pro:C_0_auto}]  We treat the
  cases (1),(2); (3); and (4),(5).

  (1),(2) From Proposition~\ref{pro:C_0_dim_md} it is clear that
  $\Aut(\mathcal{C}_0) = \Sym(\Omega)$, if $1=a \leq b$ or $2=a=b$.

  (3) For $2=a<b$ the explicit description of $\mathcal{C}_0$ in
  Proposition~\ref{pro:C_0_dim_md} shows that $\Aut(\mathcal{C}_0)$
  contains $C_2 \wr \Sym(b) = \Sym(b) \ltimes C_2^b$, where the action
  on $(i,j) \in \Omega$ of elements of the top group, respectively
  base group, of the wreath product is given by
  \begin{equation}\label{equ:wreath-action}
    \begin{split}
      (i,j)^{\sigma} = (i,j^\sigma) \qquad & \text{if $\sigma \in
        \Sym(b)$}, \\
      (i,j)^{\boldsymbol{\tau}} = (i^{\tau_j},j) \qquad & \text{if
        $\boldsymbol{\tau}=(\tau_1,\ldots,\tau_b) \in C_2^b$.}
    \end{split}
  \end{equation}
  It remains to show that the automorphism group is not larger than
  the wreath product.  Consider $\phi \in \Sym(\Omega) \setminus (C_2
  \wr \Sym(b))$.  The group $C_2 \wr \Sym(b)$ acts transitively on
  $\Omega$.  In fact, the top group acts as the full symmetric group
  on column-coordinates, and elements of the base group allow us to
  flip row-coordinates independently for each fixed column-coordinate;
  see \eqref{equ:wreath-action}.  Therefore, multiplying $\phi$ by a
  suitable element of $C_2 \wr \Sym(b)$, we may assume that
  $(1,1)^\phi = (1,1)$ and $(2,1)^\phi = (1,2)$.  Pictorially, $\phi$
  acts on a `generic' element $\sum_{i,j} c_{ij} \mathbf{e}_{ij} \in
  \mathcal{V}$ as follows:
  $$
  \begin{pmatrix}
    c_{11} & c_{12} & \cdots & c_{1b} \\
    c_{21} & c_{22} & \cdots & c_{2b}
  \end{pmatrix}^\phi =
  \begin{pmatrix}
    c_{11} & c_{21} & * & \cdots & * \\
    * & * & * & \cdots & *
  \end{pmatrix}.
  $$
  In particular, the image of $\mathbf{c}_1$ under $\phi$ is
  $$
  \mathbf{c}_1^\phi =
  \begin{pmatrix}
    1 & 0 & \cdots & 0 \\
    1 & 0 & \cdots & 0
  \end{pmatrix}^\phi =
  \begin{pmatrix}
    1 & 1 & 0 & \cdots & 0 \\
    0 & 0 & 0 & \cdots & 0
  \end{pmatrix}.
  $$
  Because the two entries in the first column of $\mathbf{c}_1^\phi$
  sum to $1$, but the two entries in the third column sum to $0$, we
  conclude that $\mathbf{c}_1^\phi \not \in \mathcal{C}_0$.  Thus
  $\phi \not \in \Aut(\mathcal{C}_0)$.

  (4),(5) Finally, we consider the case $2 < a \leq b$.  Clearly, any
  automorphism of $\mathcal{C}_0$ is uniquely determined by its effect
  on the elements of $\mathbf{R} \cup \mathbf{C}$.  

  First suppose that $2<a<b$, and recall the weight formula
  \eqref{equ:weight-formula} and the argument given in the proof of
  Proposition~\ref{pro:C_0_dim_md}.  From the latter we deduce that
  $\mathbf{C} = \{ \mathbf{v} \in \mathcal{C}_0 \mid \wt(\mathbf{v}) =
  a \}$ is invariant under $\Aut(\mathcal{C}_0)$, and furthermore that
  $$
  \mathbf{R} = \{ \mathbf{v} \in \mathcal{C}_0 \mid \wt(\mathbf{v}) =
  b \} \setminus \left\{ \sum\nolimits_{j \in Y} \mathbf{c}_j \mid Y
    \subseteq \{1,\ldots,b\} \right\}
  $$
  is $\Aut(\mathcal{C}_0)$-invariant.  Hence both sets $\mathbf{R}$
  and $\mathbf{C}$ are invariant under the action of
  $\Aut(\mathcal{C}_0)$.  Comparing with the action described in
  \eqref{equ:action}, this implies that $\Aut(\mathcal{C}_0) = \Sym(a)
  \times \Sym(b)$.

  Now consider the case $2<a=b$.  A similar argument as above shows
  that the union $\mathbf{R} \cup \mathbf{C} = \{ \mathbf{v} \in
  \mathcal{C}_0 \mid \wt(\mathbf{v}) = a \}$ is invariant under the
  action of $\Aut(\mathcal{C}_0)$.  Next observe that for any distinct
  $\mathbf{r}_i, \mathbf{r}_k \in \mathbf{R}$ and any distinct
  $\mathbf{c}_j , \mathbf{c}_l \in \mathbf{C}$ we have
  $$
  \wt(\mathbf{r}_i + \mathbf{r}_k) = \wt(\mathbf{c}_j + \mathbf{c}_l)
  = 2a, \quad \text{but } \wt(\mathbf{r}_i + \mathbf{c}_j) = 2a - 2.
  $$
  This implies that $\mathbf{R}, \mathbf{C}$ form a system of
  imprimitivity for the action of $\Aut(\mathcal{C}_0)$ on $\mathbf{R}
  \cup \mathbf{C}$.  Comparing with the action described in
  \eqref{equ:action} and noticing that ordinary matrix transposition
  yields an involution which swaps elementary row and column matrices,
  this implies that $\Aut(\mathcal{C}_0) = \Sym(a) \wr C_2$.
\end{proof}

\begin{proof}[Proof of Corollary~\ref{cor:C_0_cyclic}]
  We use without further ado the description of $\Aut(\mathcal{C}_0)$
  provided by Proposition~\ref{pro:C_0_auto}. Clearly, it suffices to
  examine the situation where $2 \leq a \leq b$, but $b \not = 2$.

  First consider the case $2=a<b$, where $\Aut(\mathcal{C}_0) = C_2
  \wr \Sym(b)$.  Let $\sigma := (1 \, 2 \, \ldots \, b) \in \Sym(b)$
  be a regular cycle in the top group, and let $\boldsymbol{\tau} :=
  (1,0,\ldots,0) \in C_2^b$ be an element of the base group acting
  non-trivially precisely in the first column.  Then $\sigma
  \boldsymbol{\tau}$ generates a regular cyclic subgroup.  Indeed, one
  checks easily that the action of $\sigma \boldsymbol{\tau}$ on
  $\Omega$ is given by $(i,j)^{\sigma \boldsymbol{\tau}} = (i,j+1)$ if
  $(i,j) \in \Omega$ with $j<b$, and by $(1,b)^{\sigma
    \boldsymbol{\tau}} = (2,1)$, $(2,b)^{\sigma \boldsymbol{\tau}} =
  (1,1)$ in the remaining two cases.

  Next consider the case $2 < a < b$, where $\Aut(\mathcal{C}_0) =
  \Sym(a) \times \Sym(b)$.  Let $(\sigma,\tau) \in \Sym(a) \times
  \Sym(b)$.  The number of orbits of $(\sigma,\tau)$ on $\Omega$ is at
  least as large as the number of orbits of $\sigma$ times the number
  of orbits of $\tau$.  Hence, if $(\sigma,\tau)$ is to generate a
  regular cyclic permutation group on $\Omega$, then $\sigma$ and
  $\tau$ are necessarily regular cycles of length $a$ and $b$,
  respectively.  But in this case the order of $(\sigma,\tau)$ is $ab
  = \lvert \Omega \rvert$ precisely if $\gcd(a,b)=1$.

  Finally consider the case $2<a=b$ where $\Aut(\mathcal{C}_0) =
  \Sym(a) \wr C_2$.  Assume for a contradiction that $\Sym(a) \wr C_2$
  contains a regular cyclic subgroup.  Then $\Sym(a) \times \Sym(a)$
  contains a cyclic subgroup $\langle (\sigma,\tau) \rangle$ with two
  orbits, each of length $a^2/2$.  Since the number of orbits of
  $(\sigma,\tau)$ is at least as large as the number of orbits of
  $\sigma$ times the number of orbits of $\tau$, we may assume without
  loss of generality that $\sigma$ is a regular cycle of length $a$
  and that $\tau$ is the product of two disjoint cycles of length
  $a-c$ and $c$, say.  Since the number of orbits of $(\sigma,\tau)$
  is two, we deduce that $1 = \gcd(a,a-c) = \gcd(a,c) = \gcd(a-c,c)$.
  Hence the order of $\tau$ is $(a-c)c$, and the order of
  $(\sigma,\tau)$ equals $(a-c)ac$.  Comparing with the orbit lengths,
  this gives $a^2/2 = (a-c)ac$.  Then $\gcd(a,a-c)=\gcd(a,c)=1$
  implies $(a,c) = (2,1)$, in contradiction to $2 < a$.
\end{proof}


\section{Binary codes associated to rectangular lattice graphs}
\label{sec:C_1}

Let $a,b \in \N$ with $a \leq b$.  We make free use of the notation
introduced in Sections~\ref{sec:prel_general}, \ref{sec:prel_C0_C1}
and \ref{sec:weight-formula}.  The aim of this section is to establish
the following results concerning the binary linear code
$\mathcal{C}_1(a,b)$ and its automorphism group.

\begin{pro}\label{pro:C_1_equal_C_0}
  If $a+b \equiv_2 1$, then $\mathcal{C}_1(a,b) = \mathcal{C}_0(a,b)$.
\end{pro}

\begin{pro}\label{pro:C_1_dim_md}
  Suppose that $a+b \equiv_2 0$.  Then $\mathcal{C}_1 =
  \mathcal{C}_1(a,b)$ has dimension $\dim(\mathcal{C}_1) =
  \dim(\mathcal{C}_0(a,b)) - 1 = a+b-2$.  Moreover, the minimum
  distance of $\mathcal{C}_1$ is $d(\mathcal{C}_1) = 2a$ if $1 \leq
  a<b$, and $d(\mathcal{C}_1) = 2a-2$ if $1<a=b$.

  In special cases we have the following explicit
  description of $\mathcal{C}_1$.
  \begin{enumerate}
  \item If $1=a=b$, then $\mathcal{C}_1 = \{0\}$.
  \item If $2=a=b$, then $\mathcal{C}_1 = \left\{ \left(
        \begin{smallmatrix} 0 & 0 \\ 0 & 0 \end{smallmatrix} \right),
      \left( \begin{smallmatrix} 0 & 1 \\ 1 & 0 \end{smallmatrix}
      \right), \left( \begin{smallmatrix} 1 & 0 \\ 0 & 1
        \end{smallmatrix} \right), \left( \begin{smallmatrix} 1 & 1 \\
          1 & 1 \end{smallmatrix} \right) \right\} $.
  \item If $1 \equiv_2 a \equiv_2 b$, then $\mathcal{C}_1 = \{
    \mathbf{v} \in \mathcal{C}_0 \mid \wt(\mathbf{v}) \equiv_2 0 \}$.
  \end{enumerate}
\end{pro}

\begin{pro}\label{pro:C_1_auto}
  Suppose that $a+b \equiv_2 0$. Let $\mathcal{C}_1 =
  \mathcal{C}_1(a,b)$ and $\mathcal{C}_0 = \mathcal{C}_0(a,b)$ as
  above.  

  Then $\Aut(\mathcal{C}_1) \leq \Aut(\mathcal{C}_0)$, with equality
  if $1 \equiv_2 a \equiv_2 b$ or $2 < a$.  In the remaining cases we
  have
  \begin{enumerate}
  \item If $2=a=b$, then $\Aut(\mathcal{C}_1) = D_8$ is a dihedral
    group of order $8$ in its natural action of degree $4$.
  \item If $2 = a < b$ and $b \equiv_2 0$, then $\Aut(\mathcal{C}_1) =
    \Sym(b) \ltimes B$, where $B = \{ (\tau_1,\ldots,\tau_b) \in C_2^b
    \mid \sum_{j=1}^b \tau_j = 0 \}$; consequently,
    $\Aut(\mathcal{C}_1)$ has index $2$ in $\Aut(\mathcal{C}_0) = C_2
    \wr \Sym(b)$.
  \end{enumerate}
\end{pro}

We recall that the structure of $\mathcal{C}_0(a,b)$ and its
automorphism group are described in Section~\ref{sec:C_0}.  We now
give the proofs of the stated results.  As before we write
$\mathcal{C}_1 = \mathcal{C}_1(a,b)$ and $\mathcal{C}_0 =
\mathcal{C}_0(a,b)$. Moreover, we define 
$$
\mathbf{a} := \sum_{i=1}^a \mathbf{r}_i = \sum_{j=1}^b \mathbf{c}_j =
  \begin{pmatrix}
 1 & 1 & \cdots & 1 \\
 \vdots & \vdots & & \vdots \\
 1 & 1 & \cdots & 1 
\end{pmatrix}.
$$

\begin{proof}[Proof of Proposition~\ref{pro:C_1_equal_C_0}]
  Suppose that $a \equiv_2 1$ and $b \equiv_2 0$.  Since the
  underlying field has characteristic $2$, we have
  $$
  \mathbf{c}_j = \mathbf{a} + (\mathbf{a} + \mathbf{c}_j) =
  \sum_{l=1}^b (\mathbf{r}_1 + \mathbf{c}_l) + \sum_{i=1}^a
  (\mathbf{r}_i + \mathbf{c}_j) \in \mathcal{C}_1 \quad \text{for $1
    \leq j \leq b$.}
  $$
  Moreover, from $\mathbf{c}_1 \in \mathcal{C}_1$ we deduce that
  $\mathbf{r}_i = (\mathbf{r}_i + \mathbf{c}_1) + \mathbf{c}_1 \in
  \mathcal{C}_1$ for $1 \leq i \leq a$.  It follows that $\mathbf{R}
  \cup \mathbf{C} \subseteq \mathcal{C}_1$, hence $\mathcal{C}_1 =
  \vspan \langle \mathbf{R} \cup \mathbf{C} \rangle = \mathcal{C}_0$.
  The argument for $a \equiv_2 0$, $b \equiv_2 1$ is very similar.
\end{proof}

\begin{proof}[Proof of Proposition~\ref{pro:C_1_dim_md}]
  The special cases $1=a=b$ and $2=a=b$ are easily dealt with.  Now
  suppose that $2<b$, and assume for the moment that we can prove the
  assertion concerning $d(\mathcal{C}_1)$.  Observe that the claimed
  value for $d(\mathcal{C}_1)$ is strictly larger than
  $d(\mathcal{C}_0) = a$; cf.\ Proposition~\ref{pro:C_0_dim_md}.  On
  the other hand, we clearly have $\mathcal{C}_1 + \vspan \langle
  \mathbf{c}_1 \rangle = \mathcal{C}_0$, and it follows that
  $\dim(\mathcal{C}_1) = \dim(\mathcal{C}_0)-1 = a+b-2$.  Moreover, we
  have $\wt(\mathbf{r}_i + \mathbf{c}_j) = a+b-2 \equiv_2 0$ for all
  $(i,j) \in \Omega$, and thus $\mathcal{C}_1 \subseteq \mathcal{W}$
  where $\mathcal{W} := \{ \mathbf{v} \in \mathcal{C}_0 \mid
  \wt(\mathbf{v}) \equiv_2 0 \}$.  Note that in the special case $1
  \equiv_2 a \equiv_2 b$ the vector $\mathbf{c}_1 \in \mathcal{C}_0$
  has weight $\wt(\mathbf{c}_1) = a \equiv_2 1$ so that
  $\dim(\mathcal{W}) = \dim(\mathcal{C}_0) -1 = \dim(\mathcal{C}_1)$.
  From this we obtain $\mathcal{C}_1 = \mathcal{W}$ as wanted.

  Hence it suffices to prove that $d(\mathcal{C}_1) = 2a$ if $1 \leq
  a<b$, and $d(\mathcal{C}_1) = 2a-2$ if $1<a=b$.  As explained, we
  shall assume that $2<b$ throughout.

  Recall the weight formula \eqref{equ:weight-formula}.  As stated in
  Subsection~\ref{sec:weight-formula}, the elements of $\mathcal{C}_1$
  are of the form $\mathbf{v} = \sum_{i \in X} \mathbf{r}_i + \sum_{j
    \in Y} \mathbf{c}_j$ where $x := \lvert X \rvert$ and $y := \lvert
  Y \rvert$ satisfy the condition $x+y \equiv_2 0$.  Since $\mathbf{a}
  + \mathbf{a} = \mathbf{0}$, we deduce for any such $\mathbf{v}$ that
  $$
  \mathbf{v} = \left( \sum\nolimits_{i \in X} \mathbf{r}_i +
    \mathbf{a} \right) + \left( \mathbf{a} + \sum\nolimits_{j \in Y}
    \mathbf{c}_j \right) = \sum\nolimits_{i \not \in X} \mathbf{r}_i +
  \sum\nolimits_{j \not \in Y} \mathbf{c}_j.
  $$
  For our analysis we may therefore assume that $x \leq \lfloor a/2
  \rfloor$.  Of course, we shall also assume that $\mathbf{v} \not =
  \mathbf{0}$.  

  If $x=0$, then $y \geq 2$ and $\wt(\mathbf{v}) = ay \geq 2a$, with
  equality if and only if $y=2$.  In this case, $\mathbf{v}=
  \mathbf{c}_j + \mathbf{c}_l$ for suitable $j \not = l$.  Likewise,
  if $x \geq 2$, then 
  $$
  \wt(\mathbf{v}) = (a-x)y + (b-y)x = (a-2x)y + (bx-2a) + 2a \geq 2a,
  $$
  with equality if and only if $(a,b,x) = (4,4,2)$ or
  $(a,x,y)=(b,2,0)$.  In the latter case, $\mathbf{v}= \mathbf{r}_i +
  \mathbf{r}_k$ for suitable $i \not = k$.  Finally, if $x=1$ (and
  consequently $1<a$ and $y \equiv_2 1$), then
  $$
  \wt(\mathbf{v}) = (a-1)y + (b-y) = (a-2)y +b.
  $$  
  The last expression takes its minimum non-zero value for $y=1$: in
  this case $\mathbf{v} = \mathbf{r}_i + \mathbf{c}_j$ for suitable
  $i,j$ and
  \begin{enumerate}
  \item[(i)] $\wt(\mathbf{v}) = 2a-2 < 2a$ if $1<a=b$;
  \item[(ii)] $\wt(\mathbf{v}) = a+b-2 \geq 2a$ if $1<a<b$, with
    equality if and only if $b = a+2$.
  \end{enumerate}
  This analysis shows that $d(\mathcal{C}_1) = 2a$ if $1 \leq a < b$,
  and $d(\mathcal{C}_1) = 2a-2$ if $1<a=b$, as wanted.
\end{proof}

\begin{proof}[Proof of Proposition~\ref{pro:C_1_auto}]
  First we treat the special cases (1) and (2).  Then we deal with the
  remaining cases (3) $1 \equiv_2 a \equiv_2 b$ and (4) $2 < a \leq b$
  with $a \equiv_2 b \equiv_2 0$.

  (1) If $2=a=b$, the group $\Aut(\mathcal{C}_1)$ is easily
  calculated from the explicit description of $\mathcal{C}_1$ in
  Proposition~\ref{pro:C_1_dim_md}.

  (2) Next consider the case $2=a<b$ and $b \equiv_2 0$.  From the
  weight analysis in the second half of the proof of
  Proposition~\ref{pro:C_1_dim_md} we see that $\mathcal{C}_1$
  contains at most two kinds of elements of minimum weight $2a=4$,
  namely those of the form $\mathbf{c}_j + \mathbf{c}_l$ and possibly
  additional elements of the form $\mathbf{r}_i + \mathbf{c}_j$.
  Therefore the set
  $$
  \mathbf{C} \dotplus \mathbf{C} := \{ \mathbf{c}_j + \mathbf{c}_l
  \mid 1 \leq j,l \leq b \text{ with } j \not = l \}
  $$
  is equal to
  $$
  \{ \mathbf{v} \in \mathcal{C}_1 \mid \exists \mathbf{w} \in
  \mathcal{C}_1 : \wt(\mathbf{v}) = \wt(\mathbf{w}) = 4 \, \land \,
  \com(\mathbf{v},\mathbf{w})= 0 \}.
  $$
  This implies that $\mathbf{C} \dotplus \mathbf{C}$ is invariant
  under $\Aut(\mathcal{C}_1)$.  Furthermore, we observe that the set
  $\mathbf{C}$ can be described as
  \begin{multline*}
    \mathbf{C} = \{ \mathbf{v} \in \mathcal{V} \mid \wt(\mathbf{v}) =
    2 \text{ and } \exists \mathbf{w}_1,\mathbf{w}_2 \in \mathbf{C}
    \dotplus \mathbf{C} : \\ \com(\mathbf{w}_1,\mathbf{w}_2) =
    \com(\mathbf{v}, \mathbf{w}_1) = \com(\mathbf{v},\mathbf{w}_2) = 2
    \}.
  \end{multline*}
  Hence $\mathbf{C}$ is invariant under $\Aut(\mathcal{C}_1)$.  As
  $\mathcal{C}_1$ is complemented in $\mathcal{C}_0$ by $\vspan
  \langle \mathbf{c}_j \rangle$ for any $j \in \{1,\ldots,b\}$, this
  shows that $\Aut(\mathcal{C}_1) \leq \Aut(\mathcal{C}_0)$.  Now
  $\Aut(\mathcal{C}_0) = C_2 \wr \Sym(b)$ by
  Proposition~\ref{pro:C_0_auto}.  Note that $\Sym(b) \ltimes B$,
  where $B = \{ (\tau_1,\ldots,\tau_b) \in C_2^b \mid \sum_{1=i}^b
  \tau_i = 0 \}$, is a subgroup of index $2$ in $C_2 \wr \Sym(b)$.
  From $\mathbf{C} \dotplus \mathbf{C} \subseteq \mathcal{C}_1$ it is
  easily seen that $\Sym(b) \ltimes B \leq \Aut(\mathcal{C}_1)$, and
  in order to prove equality it suffices to exhibit a single element
  in $C_2 \wr \Sym(b)$ which does not leave $\mathcal{C}_1$ invariant.
  Take $\boldsymbol{\tau} := (1,0,\ldots,0) \in C_2^b$.  One readily
  computes that $(\mathbf{r}_1 +\mathbf{c}_2)^{\boldsymbol{\tau}} +
  (\mathbf{r}_1 + \mathbf{c}_2) = \mathbf{c}_1 \not \in
  \mathcal{C}_1$, hence $\mathcal{C}_1$ is not invariant under the
  action of $\boldsymbol{\tau}$, as wanted.

  (3) Consider the case $1 \equiv_2 a \equiv_2 b$.  From the explicit
  description of $\mathcal{C}_1$ in Proposition~\ref{pro:C_1_dim_md}
  we see that $\Aut(\mathcal{C}_0) \leq \Aut(\mathcal{C}_1)$.  It
  remains to prove the reverse inclusion.  Since $\wt(\mathbf{a}) = ab
  \equiv_2 1$, we have $\mathcal{C}_0 = \mathcal{C}_1 + \vspan \langle
  \mathbf{a} \rangle$.  Since both summands in this decomposition are
  $\Aut(\mathcal{C}_1)$-invariant, it follows that
  $\Aut(\mathcal{C}_1) \leq \Aut(\mathcal{C}_0)$.

  (4) The final case to consider is $2 < a \leq b$ and $a \equiv_2 b
  \equiv_2 0$.  From Proposition~\ref{pro:C_0_auto} it is clear that
  $\Aut(\mathcal{C}_0) \leq \Aut(\mathcal{C}_1)$.  It remains to prove
  the reverse inclusion.  We claim that $\mathbf{R} + \mathbf{C} = \{
  \mathbf{r}_i + \mathbf{c}_j \mid (i,j) \in \Omega \}$ is invariant
  under $\Aut(\mathcal{C}_1)$.  Indeed, we contend that this set is
  equal to
  \begin{multline*}
    \mathbf{S} := \{ \mathbf{v} \in \mathcal{C}_1 \mid \wt(\mathbf{v})
    = a+b-2 \text{ and } \\ \forall \mathbf{w} \in \mathcal{C}_1:
    \wt(\mathbf{w})=a+b-2 \Rightarrow \com(\mathbf{v},\mathbf{w}) \geq
    1 \}.
  \end{multline*}
  We use again the weight formula \eqref{equ:weight-formula} and an
  analysis similar to the one in the proof of
  Proposition~\ref{pro:C_1_dim_md}.  Elements of weight $a+b-2$ in
  $\mathcal{C}_1$ are obtained from solutions $(x,y)$ of the equation
  $(a-x)y + (b-y)x = a+b-2$ satisfying the extra condition $x+y
  \equiv_2 0$.  As described earlier we may assume that $x \leq a/2$.
  We consider three cases
  \begin{enumerate}
  \item[(i)] If $x=0$, then $ay=a+b-2$ has a permissible solution $y =
    1 + (b-2)/a$ if this number is an even positive integer.  A
    corresponding element of $\mathcal{C}_1$ would have the form
    $\sum_{j \in Y} \mathbf{c}_j$ where $\lvert Y \rvert = 1 +
    (b-2)/a$.  Since $2(1+(b-2)/a) \leq b$, for any such element
    $\mathbf{v}$ there would be a similar element $\mathbf{w}$ such
    that $\com(\mathbf{v},\mathbf{w})=0$.  Hence none of these
    elements would belong to the set $\mathbf{S}$.
  \item[(ii)] If $x=1$, then $(a-1)y + (b-y) = a+b-2$ only admits the
    solution $y=1$, corresponding to elements of the form
    $\mathbf{r}_i + \mathbf{c}_j$ which we want to characterise.
  \item[(iii)] If $x \geq 2$, then $(a-x)y + (b-y)x = a+b-2$ implies
    $$
    0 = (a-2x)y + b(x-1) - a + 2 \geq 0 + b - a + 2 \geq 2,
    $$  
    a contradiction.
  \end{enumerate}
  This analysis shows that $\mathbf{R} + \mathbf{C} = \mathbf{S}$ is
  invariant under $\Aut(\mathcal{C}_1)$ as claimed.  Observe that for
  all $i,k \in \{1,\ldots,a\}$ and $j,l \in \{1,\ldots,b\}$,
  $$
  \wt((\mathbf{r}_i + \mathbf{c}_j) + (\mathbf{r}_k + \mathbf{c}_l)) =
  \begin{cases}
    2a+2b-4 & \text{if $i \not = k$ and $j \not = l$,} \\
    2a & \text{if $i=k$ and $j \not = l$,} \\
    2b & \text{if $i \not = k$ and $j=l$,} \\
    0 & \text{if $i=k$ and $j=l$.}
  \end{cases}
  $$
  Since $2 < a \leq b$, we have $0 < 2a \leq 2b < 2a+2b-4$.  Suppose
  first that $a<b$.  A simple computation shows that the set
  $\mathbf{R} + \mathbf{C}$ can be partitioned uniquely into subsets
  $\mathbf{S}_1, \ldots, \mathbf{S}_a$, each of size $b$, such that
  for every $i \in \{1,\ldots,a\}$ and all $\mathbf{v},\mathbf{w} \in
  \mathbf{S}_i$ one has $\wt(\mathbf{v} + \mathbf{w}) = 2a$.
  Moreover, we can order the sets $\mathbf{S}_1, \ldots, \mathbf{S}_a$
  such that for every $i \in \{1,\ldots,a\}$ the vector $\mathbf{r}_i$
  is characterised as the unique element $\mathbf{v} \in \mathcal{V}$
  with $\wt(\mathbf{v})=b$ and $\com(\mathbf{v},\mathbf{w})=b-1$ for
  all $\mathbf{w} \in \mathbf{S}_i$.  This shows that $\mathbf{R}$ is
  $\Aut(\mathcal{C}_1)$-invariant.  If $a=b$ a similar argument proves
  that the union $\mathbf{R} \cup \mathbf{C}$ is
  $\Aut(\mathcal{C}_1)$-invariant.  Since $\mathcal{C}_0 =
  \mathcal{C}_1 + \vspan \langle \mathbf{v} \rangle$ for any
  $\mathbf{v} \in \mathbf{R} \cup \mathbf{C}$, this implies that in
  any case $\Aut(\mathcal{C}_1) \leq \Aut(\mathcal{C}_0)$, as wanted.
\end{proof}


\section{The codes $\mathcal{K}(n_1,\ldots,n_r)$ and their
  automorphism groups}\label{sec:wreath}
Let $r \in \N$, and let $n_1,\ldots,n_r \in \N_{\geq 3}$ be odd.
Theorem~\ref{thm:wreath_auto} is an immediate consequence of the
following description of the automorphism group of the code
$\mathcal{K}_r = \mathcal{K}_r(n_1,\ldots,n_r)$, which was defined in
Section~\ref{sec:prel_K}.

\begin{thm}\label{thm:auto_K}  
  Let $r \in \N$, and let $n_1,\ldots,n_r \in \N_{\geq 3}$ be odd.
  Then $\mathcal{K}_r = \mathcal{K}_r(n_1,\ldots,n_r)$ satisfies
  $\Aut(\mathcal{K}_r) = \Sym(n_1) \wr \ldots \wr \Sym(n_r)$, where
  the wreath product acts imprimitively as a permutation group of
  degree $n_1 \cdots n_r$.
\end{thm}

We make free use of the notation introduced in
Sections~\ref{sec:prel_general} and \ref{sec:prel_K}.  Let $i \in
\{3,\ldots,r\}$ with $i \equiv_2 1$, and let $k \in \{1,\ldots,n_i\}$.
Since $\mathcal{K}_{i-1}^{(k)}$ is an isomorphic copy of
$\mathcal{K}_{i-1}$, we can use the decomposition $\mathcal{K}_{i-1} =
\mathcal{K}_{i-2}^{(1)} \oplus \ldots \oplus
\mathcal{K}_{i-2}^{(n_{i-1})}$ given by \eqref{equ:defK} to write
\begin{equation}\label{equ:decKlevel2}
\mathcal{K}_{i-1}^{(k)} = \mathcal{K}_{i-2}^{(k,1)} \oplus \ldots \oplus
\mathcal{K}_{i-2}^{(k,n_{i-1})}
\end{equation}
where $\mathcal{K}_{i-2}^{(k,m)}$, similarly defined as in
\eqref{equ:copy}, is an isomorphic copy of $\mathcal{K}_{i-2}$ with
$$
\supp(\mathcal{K}_{i-2}^{(k,m)}) = \{(k,m)\} \times \Omega_{i-2} \quad
\text{for $m \in \{1,\ldots,n_{i-1}\}$.}
$$
Our strategy for understanding the structure of $\Aut(\mathcal{K}_i)$
is based on the description of certain indecomposable elements in
$\mathcal{K}_i$.

\begin{lem}\label{lem:indecomp}
  Let $i \in \{3,\ldots,r\}$ with $i \equiv_2 1$.

  \noindent \textup{(1)} Let $k,l \in \{1,\ldots,n_i\}$ with $k \not =
  l$.  Then there exists an indecomposable element $\mathbf{v} \in
  \mathcal{K}_i$ such that
  \begin{itemize}
  \item[(i)] $\mathbf{v} \in \mathbf{a}_i^{(k,l)} +
    \mathcal{K}_{i-1}^{(k)} + \mathcal{K}_{i-1}^{(l)}$,
  \item[(ii)] $\wt(\mathbf{v}) = 2 \widehat d(\mathcal{K}_{i-1}) = 2 n_2
    n_4 \cdots n_{i-1}$,
  \item[(iii)] $\supp(\mathbf{v}) \cap \{ (k,m) \} \times \Omega_{i-2}
    \not = \varnothing$ for all $m \in \{1,\ldots,n_{i-1}\}$, 
  \item[(iv)] $\supp(\mathbf{v}) \cap \{ (l,m) \} \times \Omega_{i-2}
    \not = \varnothing$ for all $m \in \{1,\ldots,n_{i-1}\}$,  
  \end{itemize}

  \noindent \textup{(2)} Let $\mathbf{v} \in \mathcal{K}_i$, and let
  $\mathbf{a} := \pi(\mathbf{v})$ where $\pi: \mathcal{K}_i
  \rightarrow \mathcal{A}_i$ denotes the natural projection induced by
  the direct decomposition \eqref{equ:defK}.  Then 
  $$
  \wt(\mathbf{v}) \geq \left( \wt(\mathbf{a}) / \lvert \Omega_{i-1}
    \rvert \right) \cdot \widehat d(\mathcal{K}_{i-1}).
  $$
  In particular, if $\mathbf{a} \not = \mathbf{0}$, then
  $\wt(\mathbf{v}) \geq 2 \widehat d(\mathcal{K}_{i-1})$.

  \noindent \textup{(3)} If $\mathbf{v} \in \mathcal{K}_i$ is
  indecomposable with $\wt(\mathbf{v}) \leq 2 \widehat
  d(\mathcal{K}_{i-1})$, then
  \begin{itemize} 
  \item[(a)] $\mathbf{v} \in \mathbf{a}_i^{(k,l)} +
    \mathcal{K}_{i-1}^{(k)} + \mathcal{K}_{i-1}^{(l)}$ for suitable
    $k,l \in \{1,\ldots,n_i\}$ with $k \not = l$ and $\wt(\mathbf{v})
    = 2 \widehat d(\mathcal{K}_{i-1})$, or
  \item[(b)] $\mathbf{v} \in \mathcal{K}_{i-2}^{(k,m)}$ for suitable
    $k \in \{1,\ldots,n_i\}$ and $m \in \{1,\ldots,n_{i-1}\}$.
  \end{itemize}
\end{lem}

\begin{proof}
  Everything follows from the decompositions \eqref{equ:defK} and
  \eqref{equ:decKlevel2}, together with the fact that $2 \widehat
  d(\mathcal{K}_{i-1}) = 2 n_{i-1} \widehat d(\mathcal{K}_{i-2}) = 2 n_2
  n_4 \cdots n_{i-1}$; see Proposition~\ref{pro:basicsK}.
 \end{proof}

\begin{proof}[Proof of Theorem~\ref{thm:auto_K}]
  In Section~\ref{sec:prel_K} it was observed that the wreath product,
  with its natural imprimitive action, is contained in
  $\Aut(\mathcal{K}_r)$.  Hence, by induction, it suffices to show
  that the collection
  \begin{equation}\label{equ:system}
    \supp(\mathcal{K}_{r-1}^{(k)}) = \{k\} \times \Omega_{r-1}, \quad k
    \in \{1,\ldots,n_r\},
  \end{equation}
  constitutes a system of blocks for the action of
  $\Aut(\mathcal{K}_r)$ on $\Omega_r$.

  For $t \in \N \cup \{\infty\}$ we define an undirected graph
  $\Gamma_t(\mathcal{K}_r)$ with the following vertex set and edge
  set: $\mathbf{V}\Gamma_t(\mathcal{K}_r) := \Omega_r$, and
  $\mathbf{E}\Gamma_t(\mathcal{K}_r)$ consists of edges, joining
  $\omega_1$ and $\omega_2$ whenever there exists an indecomposable
  element $\mathbf{v} \in \mathcal{K}_r$ such that $\wt(\mathbf{v})
  \leq t$ and $\{\omega_1,\omega_2\} \subseteq \supp(\mathbf{v})$.
  Clearly, $\Aut(\mathcal{K}_r)$ preserves the graph structure of
  $\Gamma_t(\mathcal{K}_r)$.

  By a simple induction argument we deduce from
  Lemma~\ref{lem:indecomp} (1) that 
  \begin{enumerate}
  \item[(i)] $\Gamma_\infty(\mathcal{K}_r)$ is connected if $r
    \equiv_2 1$,
  \item[(ii)] $\Gamma_\infty(\mathcal{K}_r)$ has precisely $n_r$
    connected components if $r \equiv_2 0$.
  \end{enumerate}
  Moreover, in the case $r \equiv_2 0$, the vertex sets of the
  connected components of $\Gamma_\infty(\mathcal{K}_r)$ are exactly
  the sets listed in \eqref{equ:system} which thus form a system of
  blocks, as wanted.

  Now suppose that $r \equiv_2 1$ and put $t(r) := 2 n_2 n_4 \cdots
  n_{r-1}$.  Again by induction we draw from Lemma~\ref{lem:indecomp}
  the more precise conclusion that $\Gamma_{t(r)}(\mathcal{K}_r)$ is
  connected, while $\Gamma_{t(r) -1}(\mathcal{K}_r)$ falls into $n_r
  n_{r-1}$ connected components whose vertex sets are precisely the
  sets
  \begin{equation}\label{equ:subsystem}
    \supp(\mathcal{K}_{r-1}^{(k,m)}) = \{(k,m)\} \times \Omega_{r-1},
    \quad (k,m) \in \{1,\ldots,n_r\} \times \{1,\ldots,n_{r-1}\}.
  \end{equation}

  In order to proceed we define a variation of the graphs considered
  thus far.  For $\mathbf{v} \in \mathcal{K}_r$ let
  $\widetilde{\Gamma}_\mathbf{v}(\mathcal{K}_r)$ denote the graph
  which is obtained from $\Gamma_{t(r) -1}(\mathcal{K}_r)$ by adding
  possibly extra edges, connecting $\omega_1$ and $\omega_2$ whenever
  $\{ \omega_1, \omega_2 \} \subseteq \supp(\mathbf{v})$.  Let
  $\Delta_\mathbf{v} := \Delta_\mathbf{v}(\mathcal{K}_r)$ denote the
  vertex set of the connected component in
  $\widetilde{\Gamma}_\mathbf{v}(\mathcal{K}_r)$ which contains
  $\supp(\mathbf{v})$.
  
  Lemma~\ref{lem:indecomp} shows that, if $\mathbf{v} \in
  \mathcal{K}_r$ is indecomposable with $\wt(\mathbf{v}) = t(r)$, then
  either $\lvert \Delta_\mathbf{v} \rvert = 2 \lvert \Omega_{r-1}
  \rvert$ or $\lvert \Delta_\mathbf{v} \rvert = \lvert \Omega_{r-2}
  \lvert$.  Indeed, in the former case $\Delta_\mathbf{v}$ is the
  union of two distinct sets listed in \eqref{equ:system}, while in
  the latter case $\Delta_\mathbf{v}$ is one of the sets listed in
  \eqref{equ:subsystem}.  This implies that the sets in
  \eqref{equ:system} can be characterised as intersections
  $\Delta_{\mathbf{v},\mathbf{w}} := \Delta_\mathbf{v} \cap
  \Delta_\mathbf{w}$ where $\mathbf{v}, \mathbf{w} \in \mathcal{K}_r$
  are indecomposable with $\wt(\mathbf{v}) = \wt(\mathbf{w}) = t(r)$
  and $\lvert \Delta_{\mathbf{v},\mathbf{w}} \rvert = \lvert
  \Omega_{r-1} \rvert$.  From this description it follows that the
  sets in \eqref{equ:system} form a system of blocks, as wanted.
\end{proof}


\section{Symmetric, alternating and cyclic groups}
\label{sec:cyclic_alt}

Let $N \in \N$, and let $\mathcal{V} = \bigoplus_{i=1}^N \F_2
\mathbf{e}_i$ be an $\F_2$-vector space of dimension $N$, with fixed
standard basis $\mathbf{e}_i$ indexed by $i \in \{1,\ldots,N\}$.  Put
$\mathbf{a} := \sum_{i=1}^N \mathbf{e}_i$.  We refer to the following four
codes as \emph{elementary codes}:
$$
\mathcal{E}_0 := \{ \mathbf{0} \}, \quad \mathcal{E}_1 := \F_2
\mathbf{a}, \quad\mathcal{E}_2 := \{ \mathbf{v} \in \mathcal{V} \mid
\wt(\mathbf{v}) \equiv_2 0 \}, \quad \mathcal{E}_3 := \mathcal{V}.
$$
If $N = 1$, then $\mathcal{E}_0 = \mathcal{E}_2$ and $\mathcal{E}_1 =
\mathcal{E}_3$; if $N = 2$, then $\mathcal{E}_1 = \mathcal{E}_2$.
Otherwise the four codes are distinct.  But note that $\mathcal{E}_0$
and $\mathcal{E}_3$, respectively $\mathcal{E}_1$ and $\mathcal{E}_2$,
are orthogonal to one another with respect to the standard inner
product.  Clearly, the automorphism group of any of the elementary
codes is the full symmetric group $\Sym(N)$.  We record the following
observation and, for completeness, indicate a short proof; cf.\
\cite[Section~4]{KnSc80}.

\begin{pro} \label{pro:alt_elem}
  Let $\mathcal{C}$ be a binary linear code of length $N$.  If
  $\Alt(N) \leq \Aut(\mathcal{C})$ and $N \not = 2$, then
  $\mathcal{C}$ is one of the elementary codes $\mathcal{E}_0, \ldots,
  \mathcal{E}_3$.
\end{pro}

\begin{proof}
  Suppose that $\Alt(N) \leq \Aut(\mathcal{C})$ and $N \geq 3$.  We
  may assume that $\mathcal{C} \not \subseteq \mathcal{E}_1$ so that
  we find $\mathbf{v} \in \mathcal{C}$ with $0 < \wt(\mathbf{v}) < N$.
  Put $k := \wt(\mathbf{v})$.  As $\Alt(N)$ acts $k$-homogeneously, we
  may assume that $\mathbf{v} = \sum_{i=1}^k \mathbf{e}_i$ where $1
  \leq k < N$.  If $k=1$, then $\mathcal{C} = \mathcal{E}_3$.  If $k
  \geq 2$, then applying the $3$-cycle $\sigma := (1,k,k+1) \in
  \Alt(N)$, we see that $\mathbf{e}_1 + \mathbf{e}_{k+1} = \mathbf{v}
  + \mathbf{v}^\sigma \in \mathcal{C}$, hence $\mathcal{E}_2 \subseteq
  \mathcal{C}$.
\end{proof}

\begin{cor}\label{cor:noAlt}
  Let $\mathcal{C}$ be a binary linear code of length $N \geq 3$. Then
  $\Aut(\mathcal{C}) \not = \Alt(N)$.
\end{cor}

Now we prove Proposition~\ref{pro:cyclic_auto} and
Theorem~\ref{thm:alternating_auto}.

\begin{proof}[Proof of Proposition~\ref{pro:cyclic_auto}]
  Let $\mathcal{C}$ be a binary cyclic code of length $N$ such that $G
  = \Aut(\mathcal{C})$ is cyclic of odd order.  Since $G$ contains a
  regular cyclic subgroup of order $N$ and since any transitive cyclic
  subgroup of $\Sym(N)$ has order precisely $N$, we deduce that
  $\Aut(\mathcal{C}) = C_N$.  We may realise a code isomorphic to
  $\mathcal{C}$ as an ideal $\mathcal{I}$ of the finite ring
  $\mathcal{R} := \F_2[X] / (X^N -1)$, equipped with the standard
  basis $1,X,\ldots,X^{N-1}$.  The ideal $\mathcal{I}$ is principal
  and invariant under the Frobenius automorphism of $\mathcal{R}$.
  The latter induces a permutation $\pi$ of the standard basis, given
  by $X^m \mapsto X^{2m}$ where exponents are to be read as integers
  modulo $N$.  As $\pi$ fixes the basis element $1$, it can only
  belong to a regular cyclic group, if it is trivial.  Thus $N=1$ or
  $N=2$.  Since $N$ is odd, we conclude that $\mathcal{C}=\{0\}$, and
  $\Aut(\mathcal{C})$ is trivial.
\end{proof}  

\begin{proof}[Proof of Theorem~\ref{thm:alternating_auto}]
  Let $\mathcal{C}$ be a binary cyclic code of length $N$ such that
  $\Aut(\mathcal{C})$ is isomorphic to an alternating group $\Alt(n)$
  of degree $n \geq 3$.  An exact factorisation of $\Alt(n)$ consists
  of two subgroups $G,H \leq \Alt(n)$ such that $\Alt(n) = GH$ and $G
  \cap H = 1$.  Since $\Aut(\mathcal{C})$ contains a regular cyclic
  subgroup of order $N$ which is complemented by any point stabiliser,
  this provides an exact factorisation $\Alt(n) = GH$ with one of the
  groups $G,H$ cyclic of order $N$.  Exact factorisations of
  alternating groups were studied by Wiegold and
  Williamson~\cite{WiWi80}.  Adhering to the notation in
  \cite[Theorem~A]{WiWi80}, our setting allows for two possibilities.
  It could be that $G$ is cyclic of odd order $n=N$ and $H \cong
  \Alt(n-1)$, but this would contradict Corollary~\ref{cor:noAlt}.
  The only other possibility is that $n=8$, that $G \cong \AGL(3,2)$
  is an affine group and $H$ is cyclic of order $N=15$.  Noting that
  $\Alt(8) \cong \PSL(4,2)$, we observe that this group does indeed
  arise as the automorphism group of the binary Hamming code of length
  $2^4-1 = 15$.
\end{proof}


\section{Primitive permutation groups} \label{sec:primitive}

In this section we prove Theorem~\ref{thm:primitive}, using results
from the well-developed theory of permutation groups and modular
permutation representations.

Let $N \in \N$.  Let $G \leq \Sym(N)$ be the automorphism group of a
binary cyclic code $\mathcal{C}$, and suppose that $G$ is a primitive
permutation group.  Then $G$ contains a regular cyclic subgroup and
hence one of the following holds; see \cite[Theorem~3]{Jo02}.
\begin{enumerate}
\item $C_p \leq G \leq \AGL(1,p)$ where $p=N$ is prime.
\item $G = \Sym(N)$, or $G = \Alt(N)$ where $N \geq 3$ is odd.
\item $\PGL(d,q) \leq G \leq \PGamL(d,q)$ where $d \geq 2$, $q = p^k$
  is a prime power and $N = (q^d - 1)/(q-1)$.
\item $G = \PSL(2,11)$, $M_{11}$ or $M_{23}$ where $N = 11$, $11$ or
  $23$ respectively.
\end{enumerate}

\begin{proof}[Proof of Theorem~\ref{thm:primitive}]
  We consider the four cases listed above.

  (1) Suppose that $C_p \leq G \leq \AGL(1,p)$ where $p=N$ is prime.
  Here $\AGL(1,p)$ denotes the affine group of degree $1$.  If $p=2$
  then $G = \Sym(2)$ will be covered by case (2) below.  So suppose
  that $p \geq 3$.  Proposition~\ref{pro:cyclic_auto} shows that $C_p
  \lneqq G$. If $p=3$, then $G = \AGL(1,3) = \Sym(3)$ will be covered
  by case (2) below. Now suppose that $p \geq 5$.  For a contradiction
  assume that $G = \AGL(1,p)$.  From \cite[Table~1 and Lemma~2]{Mo80}
  we deduce that the underlying code $\mathcal{C}$ is elementary, and
  hence $G = \Aut(\mathcal{C}) = \Sym(p)$, a contradiction.

  (2) In Section~\ref{sec:cyclic_alt} it was shown that the symmetric
  group $\Sym(N)$ occurs as the automorphism group of the elementary
  codes.  According to Corollary~\ref{cor:noAlt} the alternating group
  $\Alt(N)$, $N \geq 3$, does not occur as the automorphism group of a
  binary linear code of length $N$.

  (3) Suppose that $\PGL(d,q) \leq G \leq \PGamL(d,q)$ where $d \geq
  2$, $q = p^k$ is a prime power and $N = (q^d - 1)/(q-1)$.  First we
  assume that $d=2$ and arrive at a contradiction.  As $\PGL(2,q)$
  acts $3$-transitively on $1$-dimensional projective space
  $\mathbb{P}^1(\F_q)$, we deduce from \cite[Table~1 and
  Lemma~2]{Mo80} that the underlying code $\mathcal{C}$ is elementary,
  and hence $G = \Aut(\mathcal{C}) = \Sym(N)$, a contradiction.  Hence
  $d \geq 3$, and, similarly, we deduce from \cite[Table~1 and
  Lemma~2]{Mo80} that $p=2$ and hence $q=2^k$.

  Let $\mathcal{V}$ denote the permutation module over $\F_2$,
  associated to the natural action of $\PGL(d,q)$ on
  $(d-1)$-dimensional projective space $\mathbb{P}^{d-1}(\F_q)$.  Let
  $\mathcal{U}_1$ be a $\PGL(d,q)$-submodule of $\mathcal{V}$.  We claim that
  $\mathcal{U}_1$ is automatically $\PGamL(d,q)$-invariant.  Indeed,
  let $\sigma$ be a generator of the cyclic group
  $\PGamL(d,q)/\PGL(d,q) \cong \Aut(\F_q \vert \F_2)$.  Then
  $\mathcal{U}_2 := \mathcal{U}_1^\sigma$, regarded as a
  $\PGL(d,q)$-module, is simply a twist of $\mathcal{U}_1$.  Writing
  $\overline{\F_2}$ for the algebraic closure of $\F_2$, we conclude
  that the composition factors of the $\overline{\F_2}
  \PGL(d,q)$-modules $\overline{\mathcal{U}_1} := \overline{\F_2}
  \otimes \mathcal{U}_1$ and $\overline{\mathcal{U}_2} :=
  \overline{\F_2} \otimes \mathcal{U}_2$ are the same.  The submodules
  of the $\overline{\F_2} \PGL(d,q)$-module $\overline{\mathcal{V}} :=
  \overline{\F_2} \otimes \mathcal{V}$ are uniquely determined by
  their composition factors; see \cite{BaSi00}.  Hence we conclude
  that $\overline{\mathcal{U}_1} = \overline{\mathcal{U}_2}$ and this
  implies $\mathcal{U}_1 = \mathcal{U}_2$.  We obtain $G =
  \PGamL(d,q)$, as wanted.

  For $d \geq 3$ and $q = 2^k$ we still need to justify that the
  permutation group $\PGamL(d,q)$ does indeed occur as the
  automorphism group of a suitable binary cyclic code $\mathcal{K}$.
  The explicit description of permutation modules of $\PGL(d,q)$ in
  \cite{BaSi00} guarantees the existence of a non-elementary binary
  cyclic code $\mathcal{K}$ of length $N = (q^n -1)/(q-1)$ such that
  $\PGL(d,q) \leq \Aut(\mathcal{K})$.  From \cite{KaMc74} we conclude
  that either $\Aut(\mathcal{K}) \subseteq \PGamL(d,q)$ or $\Alt(N)
  \subseteq \Aut(\mathcal{K})$.  Proposition~\ref{pro:alt_elem} rules
  out the second possibility, and hence our argument above implies
  that $\Aut(\mathcal{K}) = \PGamL(d,q)$.

  (4) A finite computation shows that of the three possible
  permutation groups precisely one, namely the Mathieu group $M_{23}$
  acting as a permutation group of degree $23$, occurs as the
  automorphism group of a binary cyclic code, namely the binary Golay
  code.
\end{proof}


\section{Examples} \label{sec:examples}

In this section we give several explicit examples of binary cyclic
codes with automorphism groups whose appearance can not be fully
explained by the results in this paper.  The examples are based on
computer calculations carried out by the first author as part of his
PhD project~\cite{Bi07}.

\subsection{Other types of groups} In~\cite[Appendices B-I]{Bi07} one
finds, in particular, a systematic listing of the non-soluble groups
which occur as automorphism groups of binary cyclic codes up to length
$70$.  In the following table we select nine single examples, in order
to illustrate that groups more complicated than those covered by
Theorems~\ref{thm:wreath_auto}, \ref{thm:direct} and
\ref{thm:primitive} occur.

\begin{table}[H]
  \begin{center}
    \begin{tabular}[h]{|c|c|l|l|}
      \hline 
      Ref.\ no.\ & $[N,k,d]$-code $\mathcal{C}$ &
      $\Aut(\mathcal{C})$ & Ref.\ in \cite{Bi07} \\ 
      \hline & & & \\[-1.2em] \hline  
      1  & $[40,35,2]$ & $\Sym(5) \wr \Sym(2) \wr \Sym(4)$ & E.11-4 \\ 
      \hline  
      2 & $[48,43,2]$ & $\Sym(6) \wr \Sym(2) \wr \Sym(4)$ & E.17-6 \\ 
      \hline  
      3  & $[48,34,2]$ & $\Sym(2) \wr \Sym(2) \wr \Sym(6) \wr \Sym(2)$
      & E.17-17 \\  
      \hline  
      4 & $[35,20,3]$ & $\PSL(3,2) \wr \Sym(5)$ & E.7-6\\
      \hline  
      5  & $[35,28,4]$ & $\Sym(5) \times \PSL(3,2)$ & E.7-4 \\
      \hline  
      6  & $[45,27,4]$ & $(\Sym(5) \times \Sym(3)) \wr \Sym(3)$ & E.15-16\\ 
      \hline  
      7  & $[48,32,4]$ & $(\Sym(8) \wr \Sym(2)) \times \Sym(3)$ & E.17-21\\ 
      \hline  
      8  & $[35,20,6]$ & $\Sym(5) \times (C_3 \ltimes C_7)$ & E.7-7 \\
      \hline  
      9  & $[85,73,4]$ & $\Sym(5) \times (C_8 \ltimes C_{17})$ & cf.\ B.3 \\
    \hline
    \end{tabular}
  \end{center}
  \caption{Selected examples of automorphism groups of binary cyclic codes}
\end{table}

The first three examples indicate that one should be able to
generalise Theorem~\ref{thm:wreath_auto} to include symmetric groups
of even degree.  However, the experimental evidence in \cite{Bi07}
also suggests that the automorphism group of a binary cyclic group is
never an iterated wreath product of symmetric groups ending in
$\Sym(2) \wr \Sym(2)$.

\subsection{Automorphism groups of affine type}
Theorem~\ref{thm:primitive} gives a description of the primitive
permutation groups which occur as automorphism groups of binary cyclic
codes.  However, it remains an open problem to find out precisely
which subgroups of the affine groups $\AGL(1,p)$ occur, where $p \geq
5$ is prime.

\begin{table}[H] \label{tab:affine}
  \begin{center}
    \begin{tabular}[h]{|c|l|l|l|}
      \hline 
      $p$ & $p-1$ & param.\ $m$ such that & $[N,k,d]$-parameters
      for codes $\mathcal{C}$ \\
      & & $\Aut(\mathcal{C}) \cong C_m \ltimes C_p$ & (in pairs,
      corresponding to dual codes) \\ 
      \hline & & & \\[-1.2em] \hline  

      $17$ & $16 = 2^4$ & $8 = 2^3$          & \hspace*{4mm}$[17,8,6]$,
      \hspace*{4mm}$[17,9,5]$ \\
      \hline  

      $31$ & $30 = 2\!\cdot\! 3 \!\cdot\! 5$ & $5$ & $[31,10,12]$,
      \hspace*{2mm}$[31,21,5]$; $[31,11,11]$, \hspace*{2mm}$[31,20,6]$; \\
      & & &  \hspace*{2mm}$[31,15,6]$, \hspace*{2mm}$[31,16,5]$;
      \hspace*{2mm}$[31,15,8]$, \hspace*{2mm}$[31,16,6]$ \\ 
      \cline{3-4}  &  & $10 = 2 \!\cdot\! 5$ & $[31,10,10]$,
      \hspace*{2mm}$[31,21,5]$; $[31,11,10]$, \hspace*{2mm}$[31,20,6]$ \\
      \cline{3-4}  &  & $15 = 3 \!\cdot\! 5$ & \hspace*{2mm}$[31,15,8]$,
      \hspace*{2mm}$[31,16,7]$ \\
      \hline  

      $41$ & $40 = 2^3\!\cdot\! 5$ & $20 = 2^2 \!\cdot\! 5$ &
      $[41,20,10]$, \hspace*{2mm}$[41,21,9]$ \\
      \hline

      $43$ & $42 = 2\!\cdot\! 3 \!\cdot\! 7$ & $14 = 2 \!\cdot\! 7$ &
      $[43,14,14]$, \hspace*{2mm}$[43,29,6]$; $[43,15,13]$,
      \hspace*{2mm}$[43,28,6]$  \\  
      \hline  

      $47$ & $46 = 2\!\cdot\! 23 $ & $23$ & $[47,23,12]$, $[47,24,11]$ \\
      \hline  

      $71$ & $70 = 2\!\cdot\! 5 \!\cdot\! 7$ & $35=5 \!\cdot\! 7$ &
      $[71,35,12]$, $[71,36,11]$ \\ 
      \hline  

      $73$ & $72 = 2^3\!\cdot\! 3^2$ & $9 = 3^2$ & \hspace*{2mm}$[73,9,28]$,
      \hspace*{2mm}$[73,64,3]$; $[73,10,28]$, \hspace*{2mm}$[73,63,4]$; \\
      & & & $[73,18,24]$, \hspace*{2mm}$[73,55,6]$; 
      $[73,19,21]$, \hspace*{2mm}$[73,54,6]$; \\
      & & & $[73,27,16]$, \hspace*{2mm}$[73,46,8]$; $[73,27,16]$,
      \hspace*{2mm}$[73,46,9]$; \\
      & & & $[73,27,18]$, \hspace*{2mm}$[73,46,9]$; $[73,27,20]$,
      \hspace*{2mm}$[73,46,9]$; \\  
      & & &  $[73,28,13]$, \hspace*{2mm}$[73,45,8]$; $[73,28,16]$,
      $[73,45,10]$; \\
      & & & $[73,28,17]$, $[73,45,10]$; $[73,36,10]$,
      \hspace*{2mm}$[73,37,9]$; \\ 
      & & &  $[73,36,12]$,
      \hspace*{2mm}$[73,37,9]$; $[73,36,12]$, $[73,37,10]$; \\ 
      & & &  $[73,36,14]$,
      \hspace*{2mm}$[73,37,9]$; $[73,36,14]$, $[73,37,12]$; \\ 
      & & & $[73,36,14]$, $[73,37,13]$ \\
      \cline{3-4}   &  & $18 = 2 \!\cdot\! 3^2$ &    $[73,18,24]$,
      \hspace*{2mm}$[73,55,6]$; $[73,19,19]$, $[73,54,6]$; \\
      &  &   &                     $[73,36,12]$, $[73,37,12]$ \\ 
      \cline{3-4}   &  & $36 = 2^2 \!\cdot\! 3^2$ & $[73,36,14]$,
      $[73,37,13]$  \\  
      \hline  

      $79$ & $78 = 2 \!\cdot\! 3 \!\cdot\! 13$ & $39 = 3 \!\cdot\! 13$ &
      $[79,39,16]$, $[79,40,15]$ \\
      \hline
    \end{tabular}
  \end{center}
  \caption{Subgroups of $\AGL(1,p)$ as automorphism groups of binary
    cyclic codes of length $N = p$ in the range $5 \leq p \leq 79$}
\end{table}

Let $p$ be a prime.  Then every binary cyclic code of length $N=p$ can
be realised as an ideal of the residue class ring $R :=
\mathbb{F}_2[X]/(X^p-1)$.  Let $f$ denote the order of $2$ in the
multiplicative group $\mathbb{F}_p^*$ and put $e := (p-1)/f$.  Then
$X^p-1$ factorises over $\mathbb{F}_2$ as a product of $X-1$ and $e$
distinct irreducible polynomials of degree $f$.  Accordingly, $R$
decomposes as a direct sum of the field $\mathbb{F}_2$ and $e$ copies
of the field $\mathbb{F}_{2^f}$.  Hence the number of ideals of $R$ is
$2^{e+1}$.  The binary cyclic codes corresponding to these ideals fall
into a certain number of isomorphism classes.  The dimensions of the
codes range over the values $kf$ and $kf+1$, where $k \in
\{0,\ldots,e\}$.  In the special case where $e=1$, there are only four
codes, namely the elementary codes discussed in
Section~\ref{sec:cyclic_alt}.  Considering a different example, if $p
= 2^l-1$ is a Mersenne prime, then $e=(p-1)/l$ and $f=l$ so that there
are $2^{((2^l-2)/l) + 1} \approx 2^{p/\log(p)}$ ideals and a priori an equal
number of corresponding binary cyclic codes of length $p$ to consider.
Clearly, as $p$ increases efficient algorithms are required to study
such a large number of codes.

Computer calculations show that, in the range $5 \leq p \leq 79$,
there exists a binary cyclic code $\mathcal{C}$ of prime length $N=p$
such that $\Aut(\mathcal{C})$ is a subgroup of the affine group
$\AGL(1,p)$ if and only if $p \in \{17,31,41,43,47,71,73,79\}$.  For
primes $p$ in this range, Table~\ref{tab:affine} lists the basic
parameters $[N,k,d]$ of all binary cyclic codes $\mathcal{C}$ of
length $N=p$ such that $\Aut(\mathcal{C}) \cong C_m \ltimes C_p$ is a
subgroup of the affine group $\AGL(1,p)$.  For convenience the prime
factorisations of $p-1$ and $m$ are exhibited.  For $p \in \{5, 7, 11,
13, 19, 23, 29, 37, 53, 59, 61, 67\}$ there exists no binary cyclic
code $\mathcal{C}$ of prime length $N = p$ such that
$\Aut(\mathcal{C})$ is a subgroup of $\AGL(1,p)$.  Except for $p = 7$
and $p = 23$, this fact can be explained by the observation that the
polynomial $X^p - 1$ admits over $\F_2$ only one irreducible factor in
addition to the trivial factor $X-1$: according to the argument given
above this implies that all binary cyclic codes of the lengths in
question are elementary.


\end{document}